\documentclass{elsart3-1}

\usepackage{graphicx}

\usepackage{amssymb}

\usepackage[english,francais]{babel}

\newtheorem{e-proposition}[theorem]{Proposition}

\newtheorem{e-definition}[theorem]{Definition\rm}

\renewcommand{\theequation}{\arabic{equation}}

\def\og{\leavevmode\raise.3ex\hbox{$\scriptscriptstyle\langle\!\langle$~}}
\def\fg{\leavevmode\raise.3ex\hbox{~$\!\scriptscriptstyle\,\rangle\!\rangle$}}

\journal{the Acad\'emie des sciences}
\begin{document}
\centerline{}
\begin{frontmatter}


\selectlanguage{english}
\title{ Relative entropy for compressible Navier-Stokes equations \\ with 
        density dependent viscosities and applications}


\selectlanguage{english}
\author[authorlabel1]{Didier Bresch\thanksref{lab1}},
\ead{Didier.Bresch@univ-savoie.fr}
\author{Pascal Noble\thanksref{lab2}},
\ead{pascal.noble@math.univ-toulouse.fr}
\author[authorlabel2]{Jean--Paul Vila}
\ead{vila@insa-toulouse.fr}
\thanks[lab1]{Research of D.B. was partially supported by the ANR project DYFICOLTI ANR-13-BS01-0003- 01}
\thanks[lab2]{Research of P.N. was partially supported by the ANR project BoND ANR-13-BS01-0009-01 .}
\address[authorlabel1]{LAMA -- UMR5127 CNRS, Bat. Le Chablais, Campus Scientifique, 73376 Le Bourget du Lac, France}
\address[authorlabel2]{IMT, INSA Toulouse, 135 avenue de Rangueil, 31077 Toulouse Cedex 9, France}


\medskip
\begin{center}
{\small Received *****; accepted after revision +++++\\
Presented by £££££}
\end{center}

\begin{abstract}
\selectlanguage{english}
   Recently A. {\sc Vasseur} and C. {\sc Yu} have proved (see \cite{VaYu}) the existence of global entropy-weak solutions to  the compressible Navier-Stokes equations with viscosities $\nu(\varrho) = \mu \varrho$ and $\lambda(\varrho) = 0$ and a pressure law under the form $p(\varrho)= a\varrho^\gamma$ with $a>0$ and $\gamma >1$ constants.
   In this note, we propose a non-trivial relative entropy for such system in a periodic box and give some applications. This extends, in some sense, results with constant viscosities recently  initiated by {\sc E.~Feireisl, B.J. Jin} and  {\sc A.~Novotny}  in  [\it J. Math. Fluid Mech. \rm (2012)].
      We present some mathematical results related to the weak-strong uniqueness,  the convergence
  to a dissipative solution of compressible or incompressible Euler equations. As a by-product,
  this mathematically justifies the convergence of solutions of a viscous shallow-water system to 
 solutions of the inviscid shallow-water system.  

\smallskip

\selectlanguage{francais}
\noindent{\bf R\'esum\'e} \vskip 0.5\baselineskip \noindent
{\bf Entropie relative pour Navier--Stokes compressible avec viscosit\'es d\'ependant de la densit\'e. }
 R\'ecem\-ment A. {\sc Vasseur} et C. {\sc Yu} ont prouv\'e (see \cite{VaYu})  l'existence globale de solutions faibles--entropiques de Navier-Stokes compressible avec  des viscosit\'es $\nu(\varrho) = \mu \varrho$, $\lambda(\varrho) = 0$  et une pression du type  $p(\varrho)= a\varrho^\gamma$ avec 
 $a>0$ et $\gamma >1$ deux constantes.
   Dans cette note, on propose une entropie relative originale pour un tel syst\`eme avec cette d\'ependance des viscosit\'es en la densit\'e et on donne quelques applications. Ceci \'etend les r\'esultats avec viscosit\'es constantes initi\'es par  {\sc E.~Feireisl, B.J. Jin} and  {\sc A.~Novotny}  dans   [\it J. Math. Fluid Mech. \rm (2012)].
 On pr\'esente quelques r\'esultats li\'es \`a l'unicit\'e faible-fort, la convergence vers une solution dissipative d'Euler compressible. 
 Ceci justifie en particulier la convergence d'un syst\`eme de Saint-venant avec viscosit\'e vers son analogue non visqueux. 
\end{abstract}
\end{frontmatter}

\selectlanguage{francais}
\section*{Version fran\c{c}aise abr\'eg\'ee}
    Dans cette note, on prolonge \`a un cas de viscosit\'es d\'ependant de la densit\'e, certains r\'esultats
maintenant connus dans le cas de viscosit\'es constantes. Plus pr\'ecisemment, on montre que toute
solution $\kappa$-entropique au sens de \cite{BrDeZa} satisfait une entropie relative qui permet
d'obtenir par exemple un r\'esultat d'unicit\'e fort-faible, un r\'esultat de convergence de 
Navier-Stokes compressible vers Euler compressible.  On am\'eliore ici l'entropie relative
introduite dans \cite{Ha} qui demandait l'hypoth\`ese forte d'une viscosit\'e proportionnelle
\`a la pression et qui consid\'erait la dimension un d'espace. Nous nous focaliserons sur le cas $\nu(\varrho)= \mu \rho$ et $\lambda(\varrho)=0$ (associ\'e aux \'equations de Saint--Venant)
pour lequel un r\'esultat complet d'existence globale de solutions $\kappa$-entropique existe: construction de solutions approch\'ees et stabilit\'e sans terme suppl\'ementaire dans le syst\`eme final.
 Plus pr\'ecisemment {\sc A. Vasseur} et {\sc C. Yu} ont obtenu (cf \cite{VaYu}) l'existence globale de
solutions faibles de Navier-Stokes compressible   satisfaisant la BD-entropie avec viscosit\'eŽ $\nu(\varrho) = \mu \varrho$ et $\lambda(\varrho)=0$: solution qui est $\kappa$-entropique pour tout
$\kappa \in (0,1)$ au sens donn\'e dans \cite{BrDeZa}. On pr\'esente quelques r\'esultats li\'es \`a l'unicit\'e faible-fort, la convergence vers une solution dissipative d'Euler compressible. Ceci justifie en particulier  le lien entre un syst\`eme de type Saint-Venant avec viscosit\'e et le syst\`eme de Saint-Venant non visqueux.  On mentionne \'egalement la convergence vers Euler incompressible. 
   La  principale difficult\'e est la d\'eg\'en\'erescence de la viscosit\'e quand la densit\'e s'annule qui ne permet pas le m\^eme contr\^ole sur la vitesse qu'avec des viscosit\'es constantes.
   La modulation de la $\kappa$-entropie introduite par le premier auteur et ses collaborateurs demande aussi un traitement pr\'ecis du terme de pression: voir ligne 3 de l'entropie relative (4) et
relation (5).
  On renvoie le lecteur interess\'e \`a \cite{BrNoVi} pour le d\'etail des preuves et une discussion sur l'extension au 
  cas o\`u les viscosit\'es $\mu$ et $\lambda$ satisfont la relation alg\'ebrique $\lambda(\varrho)=2(\mu'(\varrho)\varrho- \mu(\varrho))$: relation introduite pour la premi\`ere fois par le premier auteur et B. {\sc Desjardins} dans \cite{BrDe1}.

\selectlanguage{english}
\section{Introduction}
\label{sec1}
Since the pioneering work of {\sc C. Dafermos} and of {\sc H.-T. Yau},  relative entropy methods have become a popular and crucial tool in the study of asymptotic limits
and large-time behavior for nonlinear PDEs. 
    In a recent paper {\sc E. Feireisl, B.J. Jin, A. Novotny} (see \cite{FeJiNo}) have 
 introduced relative entropies, suitable weak solutions and weak-strong uniqueness
 for the compressible Navier-Stokes equations with constant viscosities.
   The interested reader is referred to  \cite{No}, \cite{Fe} and references
 cited therein. Based on such relative entropies, various papers have 
been dedicated to singular perturbations, see the interesting book
\cite{FeNo}, the articles \cite{BaNg}, \cite{Su} for instance.   See also 
 the recent interesting work by {\sc Th. Gallou\"et, R.~Herbin, D.~Maltese, A. Novotny}
 in \cite{GaHeMaNo} where relative entropy technics are developed to obtain 
error estimates for a numerical approximation of the compressible  Navier-Stokes 
equation with constant viscosities.
 
    Here we focus on an adaptation of the results found in \cite{FeJiNo}, \cite{BaNg} and \cite{Su} to the case of compressible Navier-Stokes equations with {\it degenerate viscosities} depending on the density
 and set in a periodic box. Relative entropy for the one-dimensional compressible Navier-Stokes equations with degenerate density dependent viscosity has been, for instance, recently  used by {\sc B. Haspot} in \cite{Ha} under the strong assumption that  the viscosity function $\nu(\varrho)$ is equal to the pressure law $p(\varrho)$ up to a multiplicative constant. This is due to the form of the modulated term in the relative entropy chosen for the  quantity coming from the pressure law.
   The main objective is to get rid of such an assumption and to extend the result to the multi-dimensional in space case. For that purpose, we will take advantage of the recent $\kappa$-entropy introduced recently by the first author, B. {\sc Desjardins} and E. {\sc Zatorska} in \cite{BrDeZa}.
    We introduce a relative entropy based on a new modulated quantity for the term involving the pressure which allows to relax the relation between the viscosity and the pressure required in \cite{Ha}. 
      By this way we cover  the (physical) case $\nu(\varrho)= \mu \varrho$, $\lambda(\varrho)=0$ and the general pressure $p(\varrho) = a\varrho^\gamma$
with $\gamma >1$.  This corresponds to the recent system for which A. {\sc Vasseur}
and {\sc C. Yu} have obtained recently global existence of entropy-weak solutions based on
the BD-entropy, the Mellet-Vasseur estimates and some original renormalization technics
linked to the momentum equations to get rid extra terms used in previous works
(drag and capillary terms). Note that  such relative entropy will be used to  prove convergence of  appropriate schemes  for the compressible  Navier--Stokes  equation with  degenerate viscosities in the forthcoming paper \cite{BrNoViVi}.
     Finally, we present some mathematical results related to the weak-strong uniqueness,  the convergence  to a dissipative solution of compressible Euler equations. 
     As a by-product,  this mathematically justifies the vanishing viscosity limit passage in viscous shallow-water equations systems.

 \section{The degenerate compressible Navier-Stokes equations and the $\kappa$-entropy}
\label{sec2}
 Let us recall the compressible Navier-Stokes equations with $\nu(\varrho)=\mu \rho$, $\lambda(\varrho)=0$:
  \begin{equation}\label{VCNS}
     \left\{  
   \begin{array}{l}
            \vspace{0.2cm}
       \partial_t \varrho + {\rm div} (\varrho{\bf u}) = 0,\\
            \vspace{0.2cm}
      \partial_t (\varrho {\bf u}) + {\rm div} (\varrho{\bf u}\otimes u)
        +\nabla p(\varrho)  -   2\mu {\rm div}(\varrho D({\bf u}))  = 0
\end{array}\right.
\end{equation}
with $D({\bf u}) = (\nabla {\bf v} + (\nabla {\bf v})^T)/2$.  Recently in \cite{BrDeZa},  it has been observed that such compressible Navier-Stokes system may be reformulated through an augmented system. 
 Introducing the intermediate velocity ${\bf v} = {\bf u} +2\kappa\mu\nabla\log(\varrho)$, a drift velocity ${\bf w} = 2\sqrt{\kappa(1-\kappa)} \mu\nabla\log(\varrho)$ and a mixture coefficient $\kappa$,  it reads
 \begin{equation}\label{SWEm}
     \left\{  
   \begin{array}{l}
            \vspace{0.2cm}
   \partial_t \varrho + {\rm div} (\varrho \left({\bf v} -2\kappa\mu \nabla\log\varrho\right))= 0,\\
            \vspace{0.2cm}
        \partial_t (\varrho {\bf v})+ {\rm div} ({\varrho {\bf v}\otimes({\bf v} -2\kappa\mu \nabla\log\varrho}) 
          +\nabla p(\varrho)= 
      \mu {\rm div}(2\varrho (1-\kappa) D({\bf v})) + \mu {\rm div}(2\kappa\varrho A({\bf v})) \\
              \vspace{0.2cm}
        \hskip7cm   -\mu {\rm div}\left(2\sqrt{\kappa(1-\kappa)}\varrho \nabla {\bf w}\right),\\
      \vspace{0.2cm}
\partial_t (\varrho{\bf w}) + {\rm div} ({\varrho {\bf w} \otimes({\bf v} -2\kappa\mu\nabla\log\varrho )})=
\mu {\rm div}(2\kappa\varrho \nabla {\bf w})
- \mu {\rm div}(2\sqrt{\kappa(1-\kappa)}\varrho  (\nabla {\bf v})^T).
\end{array}\right.
\end{equation}
with $A({\bf v}) = (\nabla {\bf v} - (\nabla {\bf v})^T)/2$. The associated $\kappa$-entropy  reads for all $t\in [0,T]$:
 \begin{equation} \label{KE-Shallow}
\begin{array}{l}
\sup_{\tau \in [0,t]} \int_\Omega\Bigl[{\varrho(\frac{|{\bf v}|^2}{2}
+ \frac{|{\bf w}|^2}{2})}+{ F(\varrho)}\Bigr](\tau){\rm d}x+2\mu \int_0^t \int_{\Omega} \varrho\left(\kappa|A({\bf v})|^2
    +|D(\sqrt{1-\kappa}{\bf v})-  \nabla(\sqrt{\kappa}{\bf w})|^2\right)  {\rm d}x\,{\rm d}s\\
\hspace{4cm} + 2\kappa\mu  \int_0^t\int_\Omega {  \frac{p'(\varrho)}{\varrho} |\nabla \varrho|^2\,{\rm d}x\,{\rm d}s} \le \int_\Omega \Bigl[{\varrho(\frac{|{\bf v}|^2}{2}
+ \frac{|{\bf w}|^2}{2}  )}+F(\varrho)\Bigr](0)
\end{array}
\end{equation}
with $s F'(s) - F(s) = p(s)$. This $\kappa$-entropy is obtained by taking the scalar product of equation satisfied by ${\bf v}$ and ${\bf w}$  respectively with ${\bf v}$ and ${\bf w}$, adding the results and  using the mass equation. Note that the more general case with $\nu(\varrho)$ arbitrary and $\lambda(\varrho) = 
2(\nu'(\varrho)\varrho - \nu(\varrho))$ for the compressible Navier-Stokes equations is covered in \cite{BrDeZa} (but with extra drag or cold pressure terms for the existence). This allows to define appropriate global solutions named $\kappa$-entropy solutions.

\section{Relative $\kappa$-entropy for compressible Navier-Stokes with degenerate viscosities}
\label{sec3}

\noindent
Let us consider the relative energy functional, denoted $E(\rho,{\bf v},{\bf w}\big\vert r,V,W)$,  defined by 
$$E(\rho,{\bf v},{\bf w}\big\vert r,V,W) = \frac{1}{2}  \int_\Omega \varrho (|{\bf w}- W|^2 + |{\bf v}-V|^2) {\rm d}x
                                      +\int_\Omega (F(\varrho)-F(r)-F'(r)(\varrho-r)){\rm d}x$$
which measures the distance between a $\kappa$-entropic weak solution $(\varrho,{\bf v},{\bf w})$ to
any smooth enough  test function $(r,V,W)$. We can prove that any weak  solution $(\rho,{\bf v},{\bf w})$ of the augmented system satisfies the following so-called  relative entropy inequality   
{\setlength\arraycolsep{1pt}
\begin{eqnarray}
\displaystyle
&& E(\rho,v,w\big\vert r,V,W) (\tau) - E(\rho,v,w\big\vert r,V,W) (0) \nonumber\\
\displaystyle
&& + 2\kappa\mu  \int_0^\tau \int_\Omega \varrho [A({\bf v}-V)|^2  
    + 2\mu\int_0^\tau \int_\Omega \varrho |D(\sqrt{(1-\kappa)}({\bf v}-V) - \sqrt\kappa ({\bf w}-W))|^2\nonumber \\
&&    + 2\kappa\mu  \int_0^\tau \int_\Omega 
     \varrho \Bigl[p'(\varrho) \nabla\log\varrho - p'(r)\nabla\log r\Bigr]\cdot
           \Bigl[ \nabla\log\varrho - \nabla\log r\Bigr] \nonumber \\
\displaystyle
&&\le \int_0^\tau \int_\Omega \varrho 
        \Bigl(( ({\bf v} - \sqrt{\frac{\kappa}{(1-\kappa)}} {\bf w})  \cdot \nabla W) \cdot (W-{\bf w}) 
      +  ( ({\bf v} - \sqrt{\frac{\kappa}{(1-\kappa)}} {\bf w}) \cdot \nabla V) \cdot (V-{\bf v}) \Bigr) \nonumber  \\
&& + \int_0^\tau\int_\Omega \varrho \Bigl( \partial_t W\cdot (W-{\bf w}) + \partial_t V \cdot (V-{\bf v})\Bigr) 
      \nonumber \\
&&
       + \int_0^\tau \int_\Omega \partial_t F'(r) (r-\varrho)
     - \int_0^\tau \int_\Omega  \nabla F'(r)\cdot  
     \Bigl[\varrho ({\bf v} - \sqrt{\frac{\kappa}{(1-\kappa)}} {\bf w})   
     -  r (V - \sqrt{\frac{\kappa}{(1-\kappa)}}W) \Bigr]
       \nonumber \\
&& + \int_0^\tau \int_\Omega(p(r)- p(\varrho)) \, {\rm div} (V- \sqrt{\frac{\kappa}{(1-\kappa)}} W) 
    \nonumber \\
&&
- \kappa  \int_0^\tau\int_\Omega p'(\varrho)\nabla\varrho \cdot [2 \mu \frac{\nabla r}{r}
     - \frac{1}{\sqrt{(1-\kappa)\kappa}} W]
     \label{RE} \\
&&+ 2\mu \int_0^\tau \int_\Omega \varrho \Bigl(D\bigl(\sqrt{(1-\kappa)}V) - \nabla(\sqrt\kappa W)\bigr)\Bigr):
                                              \Bigl(D\bigl(\sqrt{(1-\kappa)}(V-{\bf v})) - \nabla(\sqrt\kappa (W-{\bf w})\bigr)\Bigr)
      \nonumber\\
&& \displaystyle 
      + 2\kappa\mu \int_0^\tau \int_\Omega \varrho A(V): A(V-{\bf v}) 
      + 2\kappa\mu  \int_0^\tau \int_\Omega 
      \frac{\varrho}{r} p'(r) \nabla r \cdot (\frac{\nabla r}{r}- 
      \frac{\nabla\varrho}{\varrho})  \nonumber  \\
&& + 2 \sqrt{\kappa(1-\kappa)}\mu \int_0^\tau \int_\Omega \varrho
                     \Bigl[A(W):A({\bf v}-V) - A({\bf w} - W):A(V)\Bigr]  \nonumber                     
\end{eqnarray}
}
for all $\tau \in [0,T]$ and for any pair of test functions
$$ r\in {C}^1([0,T]\times \overline \Omega), \qquad  r >0, \qquad
    V,W \in {C}^1([0,T]\times \overline\Omega).$$
By using a density argument, we can of course relax the regularity on the test functions using the regularity
of the $\kappa$-entropy solutions as it was done in \cite{FeJiNo} for the constant viscosities. 
Remark that here, we do not assume $W$ to be a gradient.  The third line is also original
compared to \cite{Ha} and allows to relax the strong constraint imposed in \cite{Ha} that the viscosity is proportional
to the pressure law and covering now the physically founded case of the shallow-water equations.

\section{Some applications}
\label{sec4}
Using the relative entropy (\ref{RE}) and the identity
 {\setlength\arraycolsep{1pt}
\begin{eqnarray}
&& \varrho [p'(\varrho) \nabla\log\varrho - p'(r) \nabla\log r] \cdot [\nabla\log\varrho - \nabla\log r]
=  \varrho p'(\varrho) |\nabla\log\varrho- \nabla\log r|^2+ \\
&&\hskip2cm \nabla[p(\varrho)-p(r) - p'(r)(\varrho-r)] \cdot \nabla\log r 
 -  \bigl[\varrho(p'(\varrho) - p'(r)) - p''(r)(\varrho -r)r\bigr] \bigl|\nabla\log r|^2 ,   \nonumber
 \end{eqnarray}
}
we can justify several mathematical results. The interested reader is referred to the forthcoming paper \cite{BrNoVi} for details of the proof. 
Let us mention two of them which extends to the density dependent viscosities
 the well-known results for constant viscosities. 

\bigskip
 
\noindent {I) \it  Weak-strong uniqueness.}
Let us consider a $\kappa$--entropy solution $(\varrho,{\bf u})$
and recall ${\bf v} = {\bf u} + 2 \kappa \mu \nabla \log \varrho$ and  ${\bf w} = 2\mu \sqrt{\kappa(1-\kappa)}\nabla \log \varrho$.
Assume that $(r,W,V)$ satisfies the augmented system  with the regularity written before 
and assume that $W=2\mu \sqrt{\kappa(1-\kappa)}\nabla \log r$. Then we prove 
that  $(\varrho,{\bf v},{\bf w})= (r,V,W)$ that  means weak-strong uniqueness property: this
gives $(\varrho,{\bf u})=(r, U)$ with $U= V - \sqrt{\kappa}W/ \sqrt{(1-\kappa)}$.
More precisely we have  the following result.
\smallskip

\begin{thm} Let $\Omega$ be a periodic box. Suppose that  $p(\varrho)= a\varrho^\gamma$ with $\gamma>1$. Let $(\varrho, {\bf u})$ be a $\kappa$-entropy solution to the compressible Navier--Stokes system {\rm (\ref{VCNS})}.
 Assume that there exists a strong solution $(r,U)$ of the compressible Navier-Stokes equations
 {\rm (\ref{VCNS})} such that the terms in {\rm (\ref{RE})} are defined with $r>0$ and
 $r\in L^2(0,T;W^{1,\infty}(\Omega))\cap L^1(0,T;W^{2,\infty}(\Omega))$.
Then we have the weak-strong uniqueness result: $(\varrho,{\bf u}) = (r,U)$.
\end{thm}

\bigskip

\noindent{II) \it Inviscid limit: convergence to dissipative solution.}
Let us recall the definition of a dissipative solution of compressible Euler equations.
Such concept has been introduced by {\sc P.--L. Lions} in the incompressible setting:
see for instance \cite{Li}. The reader is referred to \cite{Fe} and \cite{BaNg} for the
extension to  the compressible framework with constant viscosities.
  Here we deal with an example of density dependent viscosities with a dissipative solution
  target. Of course the target could be the local strong solution of the compressible Euler 
  equations similarly to \cite{Su}.

\smallskip

\noindent {\bf Definition.} The pair $(\overline\varrho, \overline u)$ is a dissipative solution 
of the compressible Euler equations if and only if $(\overline\varrho,\overline u)$ satisfies
the relative energy inequality 
$$E(\overline \varrho,\overline u, 0\big\vert r,U,0)(t) \le
    E(\overline \varrho,\overline u,0\big\vert r,U,0)(0) \exp\bigl[
     \displaystyle
      c_0(r) \int_0^t \|{\rm div} U(\tau)\|_{L^\infty(\Omega)} d\tau\bigr]
  $$
 $$
     + \int_0^t \exp \bigl[
     \displaystyle
     c_0(r) \int_s^t \|{\rm div} U(\tau)\|_{L^\infty(\Omega)}\Bigr]  \int_\Omega \varrho E(r,U) 
     \cdot (U-\overline u) \, dx ds
$$
for all smooth test functions $(r,U)$ defined on $[0,T]\times \overline \Omega)$ so that $r$
 is bounded above and below away from zero and  $(r,U)$ satisfies
$$
\displaystyle
\partial_t r + {\rm div} (r U)=0, \qquad
\partial_t U + U\cdot \nabla U + \nabla F'(r) = E(r,U)
$$
\noindent
for some residual $E(r,U)$. We prove the following result.
\smallskip

\begin{thm} Let $(\varrho_\varepsilon, {\bf u}_\varepsilon)$ be any finite
$\kappa$-entropy solution to the compressible Navier-Stokes equations {\rm (\ref{VCNS})} in the
periodic setting replacing $\mu$ by $\varepsilon$. 
 Then, any weak limit $(\overline \varrho,\overline u)$ of 
$(\varrho_\varepsilon, {\bf u}_\varepsilon)$ in the sense 
{\setlength\arraycolsep{1pt}
\begin{eqnarray}
\displaystyle
&& \varrho_\varepsilon \to  \overline \varrho \hbox { weakly in  } L^\infty(0,T;L^\gamma(\Omega)), 
\nonumber\\
\displaystyle
&& \varrho_\varepsilon |{\bf v}_\varepsilon|^2
     \to  \overline \varrho \,  |\overline{u}|^2
 \hbox { weakly in  } L^\infty(0,T;L^1(\Omega)) 
 \nonumber \\
 && \varrho_\varepsilon |{\bf w}_ \varepsilon|^2
    \to 0 \hbox{  weakly in  } L^\infty(0,T;L^1(\Omega))
 \nonumber   
\end{eqnarray}
}
with ${\bf v}_\varepsilon = {\bf u}_\varepsilon + 2\varepsilon \kappa \nabla \log \varrho_\varepsilon$ and  ${\bf w}_\varepsilon = 2 \varepsilon \sqrt{\kappa(1-\kappa)} \log \varrho_\varepsilon$ as 
$\varepsilon$ tends to zero, is a dissipative solution to the compressible Euler equations.
\end{thm}

\smallskip

As a by-product this justifies the limit between a viscous shallow-water system to the inviscid 
shallow-water system. Using the relative entropy, it is 
possible to prove the convergence of the viscous shallow-water system the incompressible
 Euler equations: low Froude and inviscid limit using the mean velocity plus the oscillating
 part as target functions (see \cite{Fe} for the constant viscosities case). 



\end{document}